\documentclass[a4,12pt]{article}
\usepackage{amsmath}
\usepackage{graphics}
\usepackage{latexsym}
\usepackage{amsfonts}
\numberwithin{equation}{section}
\renewcommand{\baselinestretch}{1.40}
\title{
Uniqueness of Viscosity Solutions for Optimal Multi-Modes\\
Switching Problem with Risk of default}
\author{Brahim EL ASRI \thanks{Universit\'e Cadi Ayyad, D\'ept. de Math\'ematiques,  FSTG,
B.P. 549, Marrakech, 40.000, Maroc. e-mails:
b.elasri@uca.ma }\,\,\, \,
}
\begin{document}
\date{}
\maketitle
\newtheorem{theo}{Theorem}
\newtheorem{problem}{Problem}
\newtheorem{pro}{Proposition}
\newtheorem{cor}{Corollary}
\newtheorem{axiom}{Definition}
\newtheorem{rem}{Remark}
\newtheorem{lem}{Lemma}
\newcommand{\brm}{\begin{rem}}
\newcommand{\erm}{\end{rem}}
\newcommand{\beth}{\begin{theo}}
\newcommand{\eeth}{\end{theo}}
\newcommand{\bl}{\begin{lem}}
\newcommand{\el}{\end{lem}}
\newcommand{\bp}{\begin{pro}}
\newcommand{\ep}{\end{pro}}
\newcommand{\bcor}{\begin{cor}}
\newcommand{\ecor}{\end{cor}}
\newcommand{\be}{\begin{equation}}
\newcommand{\ee}{\end{equation}}
\newcommand{\beq}{\begin{eqnarray*}}
\newcommand{\eeq}{\end{eqnarray*}}
\newcommand{\beqa}{\begin{eqnarray}}
\newcommand{\eeqa}{\end{eqnarray}}
\newcommand{\dg}{\displaystyle \delta}
\newcommand{\cm}{\cal M}
\newcommand{\cF}{{\cal F}}
\newcommand{\cR}{{\cal R}}
\newcommand{\bF}{{\bf F}}
\newcommand{\tg}{\displaystyle \theta}
\newcommand{\w}{\displaystyle \omega}
\newcommand{\W}{\displaystyle \Omega}
\newcommand{\vp}{\displaystyle \varphi}
\newcommand{\ig}[2]{\displaystyle \int_{#1}^{#2}}
\newcommand{\integ}[2]{\displaystyle \int_{#1}^{#2}}
\newcommand{\produit}[2]{\displaystyle \prod_{#1}^{#2}}
\newcommand{\somme}[2]{\displaystyle \sum_{#1}^{#2}}
\newlength{\inter}
\setlength{\inter}{\baselineskip} \setlength{\baselineskip}{7mm}
\newcommand{\no}{\noindent}
\newcommand{\rw}{\rightarrow}
\def \ind{1\!\!1}
\def \R{I\!\!R}
\def \cadlag {{c\`adl\`ag}~}
\def \esssup {\mbox{ess sup}}
\begin{abstract}
In this paper we study the optimal $m$-states switching problem in
finite horizon as well as infinite horizon with risk of default. We
allow the switching cost functionals and cost of default to be of
polynomial growth and arbitrary. We show uniqueness of a solution
for a system of $m$ variational partial differential inequalities
with inter-connected obstacles. This system is the deterministic
version of the Verification Theorem of the Markovian optimal
$m$-states switching problem with risk of default. This problem is
connected with the valuation of a power plant in the energy market.

\end{abstract}

\no{\bf AMS Classification subjects}: 60G40 ; 62P20 ; 91B99 ; 91B28
; 35B37 ; 49L25.
\medskip

\no {$\bf Keywords$}: Real options; Backward stochastic differential
equations; Snell envelope; Stopping times ; Switching; Viscosity
solution of PDEs; Variational inequalities.

\section {Introduction}In this work we are concerned with the following systems of
$m$ variational partial differential inequalities with
inter-connected obstacles: \be\label{sysintro} \left\{
\begin{array}{l}
\min\{v_i(t,x)- \left(\max\limits_{j\in{\cal
I}^{-i}}(-g_{ij}(t,x)+v_j(t,x))\vee
(-F_i(t,x))\right),\\\qquad\qquad -\partial_tv_i(t,x)-{\cal
A}v_i(t,x)-
\psi_i(t,x)\}=0,\,\forall (t,x)\in [0,T]\times \R^k,\,\,i\in{\cal I}=\{1,...,m\},\\
v_i(T,x)=0.
\end{array}\right.
\ee \be\label{sysintro1} \left\{
\begin{array}{l}
\min\{v_i(x)- \left(\max\limits_{j\in{\cal
I}^{-i}}(-g_{ij}(x)+v_j(x))\vee (-F_i(x))\right), rv_i(x)-{\cal
A}v_i(x)-
\psi_i(x)\}=0,\,\\\forall x\in \R^k,\,\,i\in{\cal I}=\{1,...,m\},\\
\end{array}\right.
\ee

where $g_{ij}$, $\psi_i$ and $F_i$ are continuous functions, $\cal
A$ an infinitesimal generator associated with a diffusion process
and finally ${\cal I}^{-i}:=\{1,...,i-1,i+1,...,m\}.$

These system is the deterministic version of the Verification
Theorem of the optimal multi-modes switching problem in finite
horizon and infinite horizon with risk of default. These problems,
of real option type, can be introduced with the help of the
following example:

Assume we have a power station/plant which produces electricity and
which has several modes of production, e.g., the lower, the middle
and the intensive modes. The price of electricity in the market,
given by an adapted stochastic process $(X_t)_{t \leq T}$,
fluctuates in reaction to many factors such as demand level, weather
conditions, unexpected outages and so on. On the other hand,
electricity is non-storable, once produced it should be almost
immediately consumed. Therefore, as a consequence, the station
produces electricity in its instantaneous most profitable mode known
that when the plant is in mode $i\in {\cal I}$, the yield per unit
time is given by means of $\psi_i$, switching the plant from the
mode $i$ to the mode $j$ is not free and generates expenditures
given by $g_{ij}$ and, on the other hand, cost of default
(definitely stop the production) in mode $i\in {\cal I}$, is given
by $F_i$ and possibly by other factors in the energy market. So the
manager of the power plant faces two main issues:

$(i)$ when should he decide to switch the production from its
current mode to another one?

$(ii)$ to which mode the production has to be switched when the
decision of switching is made?
\medskip

\noindent In other words she faces the issue of finding the optimal
strategy of management of the plant. This issue is in relation with
the price of the power plant in the energy market.
\medskip

Optimal switching problems for stochastic systems were studied by
several authors (see e.g. \cite{[ADPS], [BE], [BO1],[BS], [CL],
[DP], [DH], [DHP], [DZ2], [E], [E1], [EH], [HJ], [TY], [dz]} and the
references therein). The motivations are mainly related to decision
making in the economic sphere. Several variants of the problem we
deal with here, including finite and infinite horizons, have been
considered during the recent years. In order to tackle those
problems, authors use mainly two approaches. Either a probabilistic
one \cite{[DH], [DHP],[E], [E1], [HJ]} or an approach which uses
partial differential inequalities (PDIs for short) \cite{[ADPS],
[BE],[BO1],[CL],[DZ2], [EH], [dz], [TY]}.

The PDIs approach turns out to study and to solve, in some sense,
the system of $m$ PDIs with inter-connected obstacles
(\ref{sysintro}) for finite horizon and (\ref{sysintro1}) for
infinite horizon. Recently El Asri and Hamad\`ene \cite{[EH]} and El
Asri \cite{[E]} extended the work of Djehiche et al \cite{[DHP]} in
the finite horizon and infinite horizon case but allowing general
jumps. In all these works the existence of the value functions of
optimal impulse control problem and uniqueness of viscosity solution
are obtained assuming that the switching problems without risk of
default.\\


Amongst the papers which consider the same problem as ours, and in
the framework of viscosity solutions approach, the most elaborated
works are certainly the ones by Djehiche and Hamad\`ene \cite{[DH]},
on the one hand, and by Arnarson et al. \cite{[ADPS]}, on the other
hand. In \cite{[DH]}, the authors show existence of a solution for
(\ref{sysintro}). Nevertheless the paper suffers from three facts:
$(i)$ the switching problem have only two modes ; $(ii)$ the
switching cost functions $g_{ij}$ should not depend on $x$ ; $(iii)$
the problem of switching is in finite horizon. The first issue of
\cite{[DH]} has been treated by Arnarson et al. \cite{[ADPS]} since
in their paper the authors show existence of the solution for
(\ref{sysintro}) in the case when the growth of the functions
$\psi_i$ is of arbitrary polynomial type. The second issue of has
been treated by El Asri and Hamad\'ene \cite{[EH]}, since in their
paper the authors show existence and uniqueness of the solution for
(\ref{sysintro}) in the case when $F_i=-\infty$. The third issu
\cite{[DH]}, i.e. considering the case of switching problem in
infinite horizon with risk of default, was right now, according to
our knowledge, an open problem. Note that in \cite{[ADPS]}, the
question of uniqueness is addressed but in the general case still
remains open. Therefore the main objective of our work, and this is
the novelty of the paper, is to show existence and uniqueness of a
solution in viscosity sense for the systems (\ref{sysintro}) and
(\ref{sysintro1}) when the functions $\psi_i$, $g_{ij}$ and $F_i$
are continuous depending also on $x$ and satisfy an arbitrary
polynomial growth condition. We show also that the solution is
unique in the class of continuous functions with polynomial growth.

This paper is organized as follows:

In Section 2, we formulate the problem and we give the related
definitions.  In Section 3, we introduce the optimal switching
problem in finite horizon and infinite horizon under consideration
and give its probabilistic Verification Theorem. It is expressed by
means of a Snell envelope of processes. Then we introduce the
approximating scheme which enables to construct a solution for the
Verification Theorem. Moreover we give some properties of that
solution. Section 4, is devoted to the connection between the
optimal switching problem in finite horizon, the Verification
Theorem and the system of PDIs (\ref{sysintro}). This connection is
made through backward stochastic differential equations with one
reflecting obstacle in the case when randomness comes from a
solution of a standard stochastic differential equation. We provide
existence and uniqueness of  viscosity solution of (\ref{sysintro})
in the class of continuous functions which satisfy a polynomial
growth condition. Section 5, we show that the solution of
(\ref{sysintro1}) is unique in the class of continuous functions
which satisfy a polynomial growth condition.$\Box$ 

\section{Assumptions and formulation of the problem}
\textbf{In finite horizon}

Throughout this paper $T$ (resp. $k$) is a fixed real (resp.
integer) positive constant. Let us now consider the followings:
\medskip

\indent $\bf H1$:  $b:[0,T]\times\R^k\rightarrow \R^{k}$ and
$\sigma:[0,T]\times\R^k\rightarrow \R^{k\times d}$ are two
continuous functions for which there exists a constant $C\geq 0$
such that for any $t\in [0, T]$ and $x, x'\in \R^k$\be
\label{regbs1}|b(t,x)|+ |\sigma(t,x)|\leq C(1+|x|) \quad \mbox{ and
} \quad |\sigma(t,x)-\sigma(t,x')|+|b(t,x)-b(t,x')|\leq C|x-x'|\ee

\indent $\bf H2$: for $i,j \in {\cal I}=\{1,...,m\}$,
$g_{ij}:[0,T]\times\R^k\rightarrow \R$,
$F_i:[0,T]\times\R^k\rightarrow \R^+$ and
$\psi_i:[0,T]\times\R^k\rightarrow \R$ are continuous functions and
of polynomial growth, $i.e.$, there exist some positive constants
$C$ and $\mu$ such that for each $i,j\in \cal I$: \be
\label{polycond} |\psi_i(t,x)|+|F_i(t,x)|+|g_{ij}(t,x)|\leq
C(1+|x|^\mu),\, \, \forall (t,x)\in [0,T]\times \R^k. \ee \indent
$\bf H3$: Moreover we assume that there exists a constant $\alpha
>0$ such that for any $(t,x)\in [0,T]\times \R^k$, \be
\min\{g_{ij}(t,x), i,j \in {\cal I},\quad i \neq j \} \geq \alpha.
\ee This condition means that switching from one mode to another one
is not free and costs at least  $\alpha>0$.
\medskip

We now consider the following system of $m$ variational inequalities
with inter-connected obstacles:  $\forall \,\,i\in {\cal I}$ \be
\label{sysvi1} \left\{
\begin{array}{l}
\min \{v_i(t,x)- \left(\max\limits_{j\in{\cal
I}^{-i}}(-g_{ij}(t,x)+v_j(t,x))\vee
(-F_i(t,x))\right),\\\qquad\qquad\qquad\qquad\qquad\qquad\qquad\qquad-\partial_tv_i(t,x)-
{\cal A}v_i(t,x)-\psi_i(t,x)\}=0,\\
v_i(T,x)=0,
\end{array}\right.
\ee where ${\cal I}^{-i}:={\cal I}-\{i\}$ and ${\cal A}$ is the
following infinitesimal generator:
\begin{equation}
\label{generateur} {\cal A}=\frac{1}{2}\sum_{i,j=1,k}(\sigma
\sigma^*)_{ij}(t,x)\frac{\partial^2} {\partial x_i \partial
x_j}+\sum_{i=1,k} b_i(t,x)\frac{\partial}{\partial x_i}\,;
\end{equation}
hereafter the superscript $(^*)$ stands for the transpose, $Tr$ is
the trace operator and finally $<x,y>$ is the inner product of
$x,y\in \R^k$.
\medskip

The first main objective of this paper is to focus on the uniqueness
of the solution in viscosity sense of (\ref{sysvi1}). To proceed we
will precise the notion of a viscosity solution of the system
(\ref{sysvi1}). It will be done in terms of subjets and superjets.

%
\begin{axiom}
Let $v \in C((0,T)\times \R^k)$, $(t,x)$ an element of $(0,T)\times
\R^k$ and finally $S_k$ the set of $k \times k$ symmetric matrices.
We denote by $J^{2,+} v(t,x)$ (resp. $J^{2,-} v(t,x)$), the
superjets (resp. the subjets) of $v$ at $(t,x)$, the set of triples
$(p,q,X)\in \R\times \R^k \times S_k$ such that:
$$\begin{array}{c}
v(s,y)\leq v(t,x) + p(s-t)+\langle q,y-x \rangle +\frac{1}{2}\langle
X(y-x),y-x\rangle+o(|s-t|+|y-x|^2) \\
(resp.\quad v(s,y)\geq v(t,x) + p(s-t)+\langle q,y-x \rangle
+\frac{1}{2}\langle X(y-x),y-x\rangle+o(|s-t|+|y-x|^2)). \Box
\end{array}$$
\end{axiom}
Note that if $\varphi-v$ has a local maximum (resp. minimum) at
$(t,x)$, then we obviously have:
$$\left(D_t \varphi(t,x),D_x \varphi(t,x),D^{2}_{xx}\varphi(t,x)\right)
\in J^{2,-} v (t,x) \,\,\, (\mbox{resp. } J^{2,+} v (t,x)). \Box$$

We now give the definition of a viscosity solution for the system of
PDE equations with risk of default in finite horizon (\ref{sysvi1}).
\begin{axiom}Let $(v_1,...,v_m)$ be a $m$-uplet of continuous functions defined
on $[0,T]\times \R^k$, $\R$-valued and such that
$(v_1,...,v_m)(T,x)=0$ for any $x\in \R^k$. The $m$-uplet
$(v_1,...,v_m)$ is called a viscosity supersolution (resp.
subsolution) of (\ref{sysvi1}) if for any $i\in {\cal I}$, $(t,x)\in
(0,T)\times \R^k$ and $(p,q,X)\in J^{2,-} v_i (t,x)$ (resp. $J^{2,+}
v_i (t,x)$),
$$
\begin{array}{l}
\min \left\{v_i(t,x)- \left(\max\limits_{j\in{\cal
I}^{-i}}(-g_{ij}(t,x)+v_j(t,x))\vee (-F_i(t,x))\right) \right.,\\
\qquad\qquad\qquad\qquad\qquad-p -\frac{1}{2}Tr[\sigma^{*} X \sigma]
-\langle b,q \rangle-\psi_{i}(t,x)\}\geq 0\, (resp. \leq 0).
\end{array}$$

It is called a viscosity solution it is both a viscosity subsolution
and supersolution .$\Box$
\end{axiom}
As pointed out previously we will show that system (\ref{sysvi1})
has a unique solution in viscosity sense. A particular case of this
system is the deterministic version of the optimal m-states
switching problem in finite horizon with risk of default which is
well documented e.g. in \cite{[ADPS],[DH]} and which we will
describe in the next section.

$\textbf{In infinite horizon}$

Let us now consider the followings assumption:
\medskip

\indent $\bf H4$:  $b:R^k\rightarrow \R^{k}$ and
$\sigma:\R^k\rightarrow \R^{k\times d}$ are two continuous functions
for which there exists a constant $C\geq 0$ such that for any $x,
x'\in \R^k$\be \label{regbs2}|b(x)|+ |\sigma(x)|\leq C(1+|x|) \quad
\mbox{ and } \quad |\sigma(x)-\sigma(x')|+|b(x)-b(x')|\leq
C|x-x'|\ee

$\bf H5$:  for $i,j \in {\cal I}=\{1,...,m\}$,
$g_{ij}:\R^k\rightarrow \R$ is a continuous function. Moreover we
assume that there exists a constant $\alpha >0$ such that for any
$x\in
 \R^k$, \be  \frac{1}{\alpha}\leq g_{ij}(x)\leq \alpha,\quad \forall i,j \in {\cal I},\quad i
\neq j . \ee

$\bf H6$:   for $i\in {\cal I}$ $\psi_i:\R^k\rightarrow \R$ and
$F_i:\R^k\rightarrow \R^+$ are  continuous function of polynomial
growth, $i.e.$, there exist a constant $C$ and $\mu$ such that for
each $i\in \cal I$: \be \label{polycond} |\psi_i(x)| +|F_i(x)| \leq
C(1+|x|^\mu),\, \, \forall x\in \R^k. \ee

\medskip

We now consider the following system of $m$ variational inequalities
with inter-connected obstacles:  $\forall \,\,i\in {\cal I}$ \be
\label{sysvi2}
\begin{array}{l}
\min\left\{v_i(x)- \left(\max\limits_{j\in{\cal
I}^{-i}}(-g_{ij}(x)+v_j(x))\vee (-F_i(x))\right),rv_i(x)- {\cal
A}v_i(x)-\psi_i(x)\right\}=0,
\end{array}
\ee where ${\cal I}^{-i}:={\cal I}-\{i\}$, $r$ is a positive
discount factor and ${\cal A}$ is the following infinitesimal
generator:
\begin{equation}
\label{generateur2} {\cal A}=\frac{1}{2}\sum_{i,j=1,k}(\sigma
\sigma^*)_{ij}(x)\frac{\partial^2} {\partial x_i \partial
x_j}+\sum_{i=1,k} b_i(x)\frac{\partial}{\partial x_i}\,.
\end{equation}

\medskip

The second main objective of this paper is to focus on the
uniqueness of the solution in viscosity sense of (\ref{sysvi2}). We
now give the definition of a viscosity solution of the elliptic
system with inter-connected obstacles (\ref{sysvi2}).

\begin{axiom}Let $(v_1,...,v_m)$ be a $m$-uplet of continuous real-valued functions defined on $R^k$. The $m$-uplet $(v_1,...,v_m)$ is called a
viscosity supersolution (resp. subsolution) of (\ref{sysvi1}) if for
any $i\in {\cal I}$, $x\in \R^k$ and $(q,X)\in J^{2,-} v_i (t,x)$
(resp. $J^{2,+} v_i (x)$),
$$
\begin{array}{l}
\min \left\{v_i(x)-\left(\max\limits_{j\in{\cal I}^{-i}}(-g_{ij}(x)+
v_j(x))\vee(-F_i(x))\right)\right.,\\ \qquad\qquad \qquad\qquad
\qquad\qquad rv_i(x) -\frac{1}{2}Tr[\sigma^{*} X \sigma] -\langle
b,q \rangle-\psi_{i}(x)\}\geq 0\, (resp. \leq 0).
\end{array}$$

It is called a viscosity solution if it is both a viscosity
subsolution and supersolution .$\Box$
\end{axiom}

As pointed out previously we will show that system (\ref{sysvi2})
has a unique solution in viscosity sense. This system is the
deterministic version of the optimal $m$-states switching problem in
infinite horizon with default risk which is well documented in
\cite{[DH], [DHP], [E1]} and which we will describe briefly in the
next section.
\section{The  optimal $m$-states  switching problem }
\subsection{In finite horizon with risk of default}
 Let $(\Omega, {\cal F}, P)$ be a
fixed probability space on which is defined a standard
$d$-dimensional Brownian motion $B=(B_t)_{0\leq t\leq T}$ whose
natural filtration is $(\cF_t^0:=\sigma \{B_s, s\leq t\})_{0\leq
t\leq T}$. Let $ \bF=(\cF_t)_{0\leq t\leq T}$ be the completed
filtration of $(\cF_t^0)_{0\leq t\leq T}$ with the $P$-null sets of
${\cal F}$, hence $(\cF_t)_{0\leq t\leq T}$ satisfies the usual
conditions, $i.e.$, it is right continuous and complete.
Furthermore, let:

- ${\cal P}$ be the $\sigma$-algebra on $[0,T]\times \Omega$ of
$\bF$-progressively measurable sets;

- ${\cal M}^{2,k}$ be the set of $\cal P$-measurable and
$\R^k$-valued processes $w=(w_t)_{t\leq T}$ such that
$E[\int_0^T|w_s|^2ds]<\infty$  and ${\cal S}^2$  be the set of $\cal
P$-measurable, continuous processes ${w}=({w}_t)_{t\leq T}$ such
that $E[\sup_{t\leq T}|{w}_t|^2]<\infty$;

-  for any stopping time $\tau \in [0,T]$, ${\cal T}_\tau$ denotes
the set of all stopping times $\theta$ such that $\tau \leq \theta
\leq T$.
\medskip

The problem of multiple switching can be described through an
example as follows. Assume  we have a plant which produces a
commodity, $e.g.$ a power station which produces electricity. The
production activity have $m$ modes, or "definitely
closed/defaulting" indicated by $\dagger$. 
A management strategy of the plant consists, on the one hand, of the
choice of a sequence of nondecreasing stopping times
$(\tau_n)_{n\geq1}$ $(i.e. \tau_n \leq \tau_{n+1}$ and $\tau_0 = 0)$
and the stopping time $\gamma$ where the manager decides to switch
the activity from its current mode to another one or definitely stop
the production. On the other hand, it consists of the choice of the
mode $\xi_n$, a r.v. ${\cal F}_{\tau_n}$-measurable with values in
${\cal I}$, to which the production is switched at $\tau_n$ from its
current mode. Therefore the admissible management strategies of the
plant are the pairs $(\delta,\xi):=((\tau_n)_{n\geq
1},\gamma,(\xi_n)_{n\geq 1})$ and the set of these strategies is
denoted by $\cal D$.

Let now $X:=(X_t)_{0\leq t\leq T}$ be an $\cal P$-measurable,
$\R^k$-valued continuous stochastic process which stands for the
market price of $k$ factors which determine the market price of the
commodity. On the other hand, assuming that the production activity
is in mode 1 at the initial time $t = 0$, let $(u_t)_{t\leq T}$
denote the indicator of the production activity's mode at time $t\in
[0, T]$ :
\begin{equation}
u_t=\ind_{[0,\tau_1]}(t)+\sum_{n\geq1}\xi_n
\ind_{(\tau_{n},\tau_{n+1}]}(t).
\end{equation}
Then for any $t\leq T$, the state of the whole economic system
related to the project at time $t$ is represented by the vector :
\begin{equation}
\begin{array}{ll}
(t, X_t, u_t),\quad &\mbox{if}\, \tau_n<t\leq \tau_{n+1};\\
(\gamma, X_{\gamma}),\quad &\mbox{if in mode}\,\, \dagger.
\end{array}
\end{equation}

Finally, let $\psi_i(t,X_t)$ be the instantaneous profit when the
system is in state $(t, X_t, i)$, for $i,j \in {\cal I} \quad i\neq
j$, let $g_{ij}(t,X_t)$ denote the switching cost of the production
at time $t$ from current mode $i$ to another mode $j$ and let
$F_i(\gamma,X_\gamma)$ denote the cost of default (definitely stop
the production) at time $\gamma$, when in mode $i$ and denote
$F_i(\gamma,X_\gamma)=F(\gamma,X_\gamma,u_\gamma)$ when
$u_\gamma=i$. Then if the plant is run under the strategy
$(\delta,\xi)=((\tau_n)_{n\geq 1},\gamma,(\xi_n)_{n\geq 1})$ the
expected total profit is given by:
$$\begin{array}{l} J(\delta,\xi)=E[\integ{0}{\gamma}\psi_{u_s}(s,X_s)ds
-\sum_{n\geq 1} g_{u_{\tau_{n-1}}u_{\tau_n}}(\tau_{n},X_{\tau_{n}})
\ind_{[\tau_{n}<\gamma]}-F(\gamma,X_\gamma,u_\gamma)\ind_{[\gamma<T]}].
\end{array}
$$
Therefore the problem we are interested in is to find an optimal
strategy, $i.e$, a strategy $(\delta^*,u^*)=((\tau^*_n)_{n\geq
1},\gamma^*),(\xi^*_n))$ such that $J(\delta^*,\xi^*)\ge
J(\delta,\xi)$ for any $(\delta,\xi)\in \cal D$.
\medskip

Note that in order that the quantity $J(\delta,\xi)$ makes sense we
assume throughout this paper that for any $i,j\in {\cal I}$ the
processes  $(F_i(t,X_t))_{t\leq T}$ , $(g_{ij}(t,X_t))_{t\leq T}$
(resp. $(\psi_i(t,X_t))_{t\leq T}$) belong to ${\cal S}^{2}$ (resp.
${\cal M}^{2,1}$). On the other hand there is a bijective
correspondence between the pairs $(\delta,\xi)$ and the pairs
$(\delta,u)$. Therefore throughout this paper one refers
indifferently to $(\delta,\xi)$ or $(\delta,u)$.
\medskip

The verification Theorem for the $m$-states optimal switching with
risk of default problem is the following:
\begin{theo}\label{thmverif}
\noindent Assume that there exist $m$ processes $(Y^i:=(Y^i_t)_{0\le
t\leq T}, i=1,...,m)$  of ${\cal S}^2$ such that: \be
\begin{array}{ll}
\label{eqvt} \forall t\leq T,\,\,&Y^i_t=\esssup_{\tau \geq
t}E[\int_t^\tau\psi_i(s,X_s)ds +\max\limits_{j\in {\cal
I}^{-i}}(-g_{ij}(\tau,X_\tau)+Y^j_\tau)\vee F_i(\tau,X_\tau)1_{[\tau
<T]}|\cF_t],\\ &Y^i_T=0.
\end{array}
\ee Then:
\begin{itemize}
\item[$(i)$]
$ Y^1_0=\sup \limits_{(\delta,\xi) \in {\cal D}}J(\delta,u). $
\item[$(ii)$] Define the sequence of $\bF$-stopping times $\delta=(\tau_t^i,\gamma),\,i=1,2,...,m,\, 0\leq t \leq T.$ as follows :
$$
\begin{array}{lll}
\tau^i_{t}&=&\sigma^i_t\wedge\widetilde{\sigma}_t^i\wedge T,\\

\end{array}
$$
where:
\begin{itemize}
\item[$\bullet$] The first time the activity defaults while in mode $i$ is given
by$$\sigma^i_t:=\inf\{s\geq t,\, Y^i_s=F_i(s,X_s)\}\wedge T,\quad
i=1,...,m.$$
\item[$\bullet$] The first time the activity is switched from mode $i$ to any of the other modes $j\neq i$ is given by
$$\widetilde{\sigma}^i_t:=\inf\{s\geq t,\, Y^i_s=\max\limits_{j\neq
i}(-g_{ij}(s,X_s)+Y_s^j)\}\wedge T$$

\end{itemize}
Finally, let $\gamma=\sup\limits_{0\leq t\leq T}\tau^i_t.$ Then, the
strategy $(\delta=((\tau^i_t)_{t\geq 0},\gamma),u^*)$ is optimal.
$\Box$
\end{itemize}
\end{theo}
\medskip
\no {\it Proof}. The arguments of proof are standard, based on the
properties the Snell envelope and is proved in \cite{[DHP]}, Theorem
1. $\Box$

The issue of existence of the processes $Y^1,...,Y^m$ which satisfy
(\ref{eqvt}) is also addressed in \cite{[DHP]}. Also for $n\geq 0$
let us define the processes $(Y^{1,n},...,Y^{m,n})$ recursively as
follows: for $i\in {\cal I}$ we set,
\begin{equation}\label{y0}
Y^{i,0}_t=\mbox{ess sup}_{\tau\geq
t}E[\integ{t}{\tau}\psi_i(s,X_s)ds+F_i(\tau,X_\tau)1_{[\tau
<T]}|{\cal F}_t],\,\,0\le t\leq T,
\end{equation}
and for $n\geq 1$,
\begin{equation}
\label{eq24} Y^{i,n}_t=\mbox{ess sup}_{\tau\geq t}
E[\integ{t}{\tau}\psi_i(s,X_s)ds+\max\limits_{k\in {\cal
I}^{-i}}(-g_{ik}(\tau,X_\tau)+Y^{k,n-1}_\tau) \vee
F_i(\tau,X_\tau)1_{[\tau <T]}|{\cal F}_t],\,\,0\le t\leq T.
\end{equation}
Then the sequence of processes $((Y^{1,n},...,Y^{m,n}))_{n\geq 0}$
have the following properties:
\begin{pro} (\cite{[DHP]}, Pro.3 and Th.2)
\begin{itemize}
\item[$(i)$] for any $i\in {\cal I}$ and
$n\geq 0$, the processes $Y^{1,n},...,Y^{m,n}$ are well-posed,
continuous and belong to ${\cal S}^2$, and verify
\be\label{croi}\forall t\le T,\,\,Y^{i,n}_t\leq Y^{i,n+1}_t\leq
E[\int_t^T\{\max_{i=1,m}|\psi_i(s,X_s)|\}ds|{\cal F}_t];\ee

\item[$(ii)$] there exist $m$ processes $Y^1,...,Y^m$ of ${\cal S}^2$ such
that for any $i\in {\cal I}$:
\begin{itemize}
\item[$(a)$]  $\forall t\leq T$, $Y^i_t=\lim\limits_{n\rightarrow
\infty}\nearrow Y^{i,n}_t$ and
$$
E[\sup_{s\leq T}|Y^{i,n}_s-Y^i_s|^2]\to 0 \quad\mbox{ as }\quad n\to
+\infty $$
\item[$(b)$] $\forall t\leq T$,
\begin{equation}
\label{eq26}{Y}^{i}_t=\mbox{ess sup}_{\tau\geq
t}E[\integ{t}{\tau}\psi_i(s,X_s)ds+ \max\limits_{k\in {\cal
I}^{-i}}(-g_{ik}(\tau,X_\tau)+{Y}^{k}_\tau) \vee
F_i(\tau,X_\tau)1_{[\tau <T]}|{\cal F}_t]
\end{equation}
i.e. ${Y}^1,...,{Y}^m$ satisfy the Verification Theorem
\ref{thmverif}.$\Box$
\end{itemize}
\end{itemize}
\end{pro}

\begin{rem} \label{unic}Note that the characterization (\ref{eq26}) implies that the processes $Y^1,...,Y^m$ of ${\cal S}^2$ which satisfy the Verification
Theorem are unique.
\end{rem}

\subsection{ In infinite horizon with risk of default}
Let $(\Omega, {\cal F}, P)$ be a fixed probability space on which is
defined a standard $d$-dimensional Brownian motion $B=(B_t)_{t\geq
0}$ whose natural filtration is $(\cF_t^0:=\sigma \{B_s, s\leq
t\})_{t\geq0}$. Let $ \bF=(\cF_t)_{ t\geq0}$ be the completed
filtration of $(\cF_t^0)_{t\geq0}$ with the $P$-null sets of ${\cal
F}$, hence $(\cF_t)_{t\geq0}$ satisfies the usual conditions,
$i.e.$, it is right continuous and complete. Furthermore, let:

- ${\cal P}$ be the $\sigma$-algebra on $[0,+\infty)\times \Omega$
of $\bF$-progressively measurable sets;

- ${\cal M}^{2,k}$ be the set of $\cal P$-measurable and
$\R^k$-valued processes $w=(w_t)_{t\geq0}$ such that\\
  $E[\int_0^{+\infty}|w_s|^2ds]<\infty$  and ${\cal S}^2$  be the set
of $\cal P$-measurable, continuous processes ${w}=({w}_t)_{t\geq 0}$
such that $E[\sup_{t\geq0}|{w}_t|^2]<\infty$;

-  for any stopping time $\tau \in \R^+$, ${\cal T}_\tau$ denotes
the set of all stopping times $\theta$ such that $\tau \leq \theta;
$

-  for any stopping time $\tau$, $ {\cal F}_{\tau}$ is the
$\sigma$-algebra on $\Omega$ which contains the sets $A$ of
$\cal{F}$ such that $ A \cap \{\tau \leq t\}\in {\cal F}_t$ for
every $t\geq 0$.$\Box$
\medskip

A decision (strategy) of the problem of multiple switching, on the
one hand, consists of the choice of a sequence of nondecreasing
stopping times $(\tau_n)_{n\geq1}$ $(i.e. \tau_n \leq \tau_{n+1}$
and $\tau_0 = 0)$ and the stopping time $\gamma$ where the manager
decides to switch the activity from its current mode to another one
or definitely stop the production. On the other hand, it consists of
the choice of the mode $\xi_n$, a r.v. ${\cal
F}_{\tau_n}$-measurable with values in ${\cal I}$, to which the
production is switched at $\tau_n$ from its current mode. Therefore
the admissible management strategies of the plant are the pairs
$(\delta,\xi):=((\tau_n)_{n\geq 1},\gamma,(\xi_n)_{n\geq 1})$ and
the set of these strategies is denoted by $\cal D$.


Let now $X:=(X_t)_{t\geq0}$ be an $\cal P$-measurable, $\R^k$-valued
continuous stochastic process which stands for the market price of
$k$ factors which determine the market price of the commodity. On
the other hand, assuming that the production activity is in mode 1
at the initial time $t = 0$, let $(u_t)_{t\geq0}$ denote the
indicator of the production activity's mode at time $t\in \R^+$ :
\begin{equation}
u_t=\ind_{[0,\tau_1]}(t)+\sum_{n\geq1}\xi_n
\ind_{(\tau_{n},\tau_{n+1}]}(t).
\end{equation}
Then for any $t\geq0$, the state of the whole economic system
related to the project at time $t$ is given by the vector:
\begin{equation}
\begin{array}{ll}
(t, X_t, u_t),\quad &\mbox{if}\, \tau_n<t\leq \tau_{n+1};\\
(\gamma, X_{\gamma}),\quad &\mbox{if in mode}\,\, \dagger.
\end{array}
\end{equation}

Finally, let $\psi_i(X_t)$ be the instantaneous profit when the
system is in state $(t, X_t, i)$, for $i,j \in {\cal I} \quad i\neq
j$, let $g_{ij}(X_t)$ denote the switching cost of the production at
time $t$ from the current mode $i$ to another mode $j$ and let
$F_i(\gamma,X_\gamma)$ denote the cost of default (definitely stop
the production) at time $\gamma$, when in mode $i$ and denote
$F_i(X_\gamma)=F(X_\gamma,u_\gamma)$ when $u_\gamma=i$. When the
plant is run under the strategy $(\delta,\xi)=((\tau_n)_{n\geq
1},\gamma,(\xi_n)_{n\geq 1})$ the expected total profit is given by:
$$\begin{array}{l} J(\delta,\xi)=E[\integ{0}{\gamma}e^{-rs}\psi_{u_s}(X_s)ds
-\sum_{n\geq 1}
e^{-r\tau_n}g_{u_{\tau_{n-1}}u_{\tau_n}}(X_{\tau_{n}})\ind_{[\tau_{n}<\gamma]}-e^{-r\gamma}F(X_\gamma,u_\gamma)].
\end{array}
$$
Then the problem we are interested in is to find an optimal
strategy, $i.e$, a strategy $(\delta^*,\xi^*)$ such that
$J(\delta^*,\xi^*)\ge J(\delta,\xi)$ for any $(\delta,\xi)\in \cal
D$.
\medskip

Note that in order that the quantity $J(\delta,\xi)$ makes sense we
assume throughout this paper that for any $i\in {\cal I}$ the
processes $(e^{-rt}\psi_i(X_t))_{t\geq0}$ and
$(e^{-rt}F_i(X_t))_{t\geq 0}$ belong to ${\cal M}^{2,1}$ and ${\cal
S}^{2}$ respectively.
%

The Verification Theorem for the $m$-states optimal switching with
risk of default problem in infinite horizon is the following:
\begin{theo}.\label{thmverif}
\noindent Assume that there exist $m$ processes
$(Y^i:=(Y^i_t)_{t\geq0}, i=1,...,m)$  of ${\cal S}^2$ such that: \be
\begin{array}{l}
\label{eqvt} \forall t\geq0,\,\,e^{-rt}Y^i_t=\esssup_{\tau \geq
t}E[\int_t^\tau e^{-rs}\psi_i(X_s)ds +e^{-r\tau}\max\limits_{j\in
{\cal I}^{-i}}(-g_{ij}(X_\tau)+Y^j_\tau)\vee
e^{-r\tau}F_i(X_\tau)|\cF_t],\quad
\\\lim\limits_{t\rightarrow+\infty}(e^{-rt}Y^i_{t})=0.
\end{array}
\ee Then:
\begin{itemize}
\item[$(i)$]
$ Y^1_0=\sup \limits_{(\delta,\xi) \in {\cal D}}J(\delta,u). $
\item[$(ii)$] Define the sequence of $\bF$-stopping times $\delta=(\tau_t^i,\gamma),\,i=1,2,...,m,\, t\geq 0.$ as follows :
$$
\begin{array}{lll}
\tau^i_{t}&=&\sigma^i_t\wedge\widetilde{\sigma}_t^i,\\

\end{array}
$$
where:
\begin{itemize}
\item[$\bullet$] The first time the activity defaults while in mode $i$ is given
by$$\sigma^i_t:=\inf\{s\geq t,\, Y^i_s=F_i(X_s)\},\quad i=1,...,m.$$
\item[$\bullet$] The first time the activity is switched from mode $i$ to any of the other modes $j\neq i$ is given by
$$\widetilde{\sigma}^i_t:=\inf\{s\geq t,\, Y^i_s=\max\limits_{j\neq
i}(-g_{ij}(s,X_s)+Y_s^j)\}$$

\end{itemize}
Finally, let $\gamma=\sup\limits_{t\geq 0}\tau^i_t.$ Then, the
strategy $(\delta=((\tau^i_t)_{t\geq 0},\gamma)$  is optimal. $\Box$
\end{itemize}
\end{theo}
\medskip
\no {\it Proof}. The arguments of proof are standard, based on the
properties the Snell envelope and is proved in \cite{[E]},
Theorem.1. $\Box$

\medskip
The issue of existence of the processes $Y^1,...,Y^m$ which satisfy
(\ref{eqvt}) is also addressed in \cite{[DHP]}. For $n\geq 0$ let us
define the processes $(Y^{n,1},...,Y^{n,m})$ recursively as follows:
for $i\in {\cal I}$ we set,
\begin{equation}\label{y00}
e^{-rt}Y^{0,i}_t=\mbox{ess sup}_{\tau\geq
t}E[\integ{t}{\tau}e^{-rs}\psi_i(X_s)ds+e^{-r\tau}F_i(X_\tau)|{\cal
F}_t],\,\, t\geq0,
\end{equation}
and for $n\geq 1$,
\begin{equation}
\label{eq242} e^{-rt}Y^{n,i}_t=\mbox{ess sup}_{\tau\geq t}
E[\integ{t}{\tau}e^{-rs}\psi_i(X_s)ds+e^{-r\tau}\max\limits_{k\in
{\cal I}^{-i}}(-g_{ik}(X_\tau)+Y^{n-1,k}_\tau)\vee
e^{-r\tau}F_i(X_\tau) |{\cal F}_t],\,\, t\geq0.
\end{equation}
Then the sequence of processes $((Y^{n,1},...,Y^{n,m}))_{n\geq 0}$
have the following properties:
\begin{pro} (\cite{[DHP]}, Pro.3 and Th.2)
\begin{itemize}
\item[$(i)$] for any $i\in {\cal I}$ and
$n\geq 0$, the processes $Y^{n,1},...,Y^{n,m}$ are well-posed,
continuous and belong to ${\cal S}^2$, and verify
\be\label{croi}\forall t\geq0,\,\,e^{-rt}Y^{n,i}_t\leq
e^{-rt}Y^{n+1,i}_t\leq
E[\int_t^{+\infty}e^{-rs}\{\max_{i=1,m}|\psi_i(X_s)|\}ds|{\cal
F}_t];\ee

\item[$(ii)$] there exist $m$ processes $Y^1,...,Y^m$ of ${\cal S}^2$ such
that for any $i\in {\cal I}$:
\begin{itemize}
\item[$(a)$]  $\forall t\geq0$, $Y^i_t=\lim\limits_{n\rightarrow
\infty}\nearrow Y^{n,i}_t$ 
\item[$(b)$] $\forall t\geq0$,
\begin{equation}
\label{eq262}e^{-rt}{Y}^{i}_t=\mbox{ess sup}_{\tau\geq
t}E[\integ{t}{\tau}e^{-rs}\psi_i(X_s)ds+ e^{-r\tau}\max\limits_{k\in
{\cal I}^{-i}}(-g_{ik}(X_\tau)+{Y}^{k}_\tau)\vee
e^{-r\tau}F_i(X_\tau) |{\cal F}_t]
\end{equation}
i.e. ${Y}^1,...,{Y}^m$ satisfy the Verification Theorem
\ref{thmverif}.$\Box$

\end{itemize}
\end{itemize}
\end{pro}

\begin{rem} \label{unic}The characterization (\ref{eq262}) implies that the processes $Y^1,...,Y^m$ of ${\cal S}^2$ which satisfy the Verification
Theorem are unique.$\Box$
\end{rem}
\section{Uniqueness of the viscosity solution in finite horizon}
Let $(t,x)\in [0,T]\times \R^k$ and let $(X^{tx}_s)_{s\leq T}$ be
the solution of the following standard SDE:
\begin{equation}\label{sde}
dX^{tx}_s=b(s,X_s^{tx})ds+\sigma(s,X_s^{tx})dB_s \mbox{ for }t\leq
s\leq T\mbox{ and }X_s^{tx}=x \mbox{ for }s\leq t\end{equation}where
the functions $b$ and $\sigma$ are the ones of (\ref{regbs1}). These
properties of $\sigma$ and $b$ imply in particular that the process
$(X^{tx}_s)_{0\le s\leq T}$ solution of the standard SDE (\ref{sde})
exists and is unique, for any $t\in [0, T]$ and $x\in \R^k$.

The operator $\cal A$  that is appearing in (\ref{generateur}) is
the infinitesimal generator associated with $X^{t,x}$.
 \beth
There are deterministic functions $v^1,...,v^m$ $:[0,T]\times
\R^k\rightarrow \R$ such that:
$$\forall (t,x)\in [0,T]\times \R^k, \forall s\in [t,T],
Y_s^{i,tx}=v^i(s,X^{tx}_s), \,\,i=1,...,m.$$  Moreover the functions
$(v^1,...,v^m):[0,T]\times \R^k\rightarrow \R$ are continuous,
solution in viscosity sense of the system of variational
inequalities with inter-connected obstacles (\ref{sysvi1})and of
polynomial growth.$\Box$\eeth $Proof$: The continuity of functions
$v^1,...,v^m$ follows from the dynamic programming principle and is
proved in \cite{[EH]}.$\Box$

Now we give an equivalent of quasi-variational inequality
(\ref{sysvi1}). In this section, we consider the new function
$\Gamma_i$ given by the classical change of variable $\Gamma_i(t,x)
= \exp(t)v_i(t, x)$, for any $t\in[0,T ]$ and $x\in \R^k$. Of
course, the function $\Gamma_i$ is continuous and of polynomial
growth with respect to its arguments.\\ A second property is given
by the

\begin{pro}
$v_i$ is a viscosity solution of (\ref{sysvi1}) if and only if
$\Gamma_i$ is a viscosity solution to the following
quasi-variational inequality in $[0,T [\times \R^k$, \be
\label{sysvi10} \left\{
\begin{array}{ll}
\min \{\Gamma_i(t,x)- \left(\max\limits_{j\in{\cal
I}^{-i}}(-e^{t}g_{ij}(t,x)+\Gamma_j(t,x))\vee
(-e^{t}F_i(t,x))\right),\\
\qquad\qquad\qquad\qquad\qquad\qquad\Gamma_i(t,x)-\partial_t\Gamma_i(t,x)-
{\cal A}\Gamma_i(t,x)-e^{t}\psi_i(t,x)\}=0,\\
\Gamma_i(T,x)=e^{T}v_i(T,x)=0.\Box
\end{array}\right.
\ee
\end{pro}

We are going now to address the question of uniqueness of the
viscosity solution of the system (\ref{sysvi1}). We have the
following:

\beth \label{uni}The solution in viscosity sense of the system of
variational inequalities with inter-connected obstacles
(\ref{sysvi1}) is unique in the space of continuous functions on
$[0,T]\times R^k$ which satisfy a polynomial growth condition, i.e.,
in the space
$$\begin{array}{l}{\cal C}:=\{\varphi: [0,T]\times \R^k\rightarrow
\R, \mbox{ continuous and for any }\\\qquad \qquad\qquad(t,x), \,
|\varphi(t,x)|\leq C(1+|x|^\mu) \mbox{ for some constants } C \mbox{
and }\mu\}.\end{array}$$ \eeth {\it Proof}. We will show by
contradiction that if $u_1,...,u_m$ and $w_1,...,w_m$ are a
subsolution and a supersolution respectively for (\ref{sysvi10})
then for any $i=1,...,m$, $u_i\leq w_i$. Therefore if we have two
solutions of (\ref{sysvi10}) then they are obviously equal. Actually
for some $R>0$ suppose there exists
$(\overline{t},\overline{x},\overline{i})\in(0,T)\times B_R\times
{\cal I}$ $(B_R := \{x\in \R^k; |x|<R\})$ such that:
\begin{equation}
\label{comp-equ}
\max\limits_{t,x,i}(u_i(t,x)-w_i(t,x))=u_{\overline{i}}(\overline{t},\overline{x})-w_{\overline{i}}(\overline{t},\overline{x})=\eta>0.
\end{equation}Let us take
$\theta,\lambda$ and $\beta \in (0,1]$ small enough, so that the
following holds: \be \left\{
\begin{array}{llll}
\beta T<\frac{\eta}{4}\\
-\lambda w_{\overline{i}}(\overline{t},\overline{x})< \frac{\eta}{4}\\
\frac{\lambda}{\overline{t}}<\frac{\eta}{4}.
\end{array}
\right. \ee Here $\gamma$ is the growth exponent of the functions
which w.l.o.g we assume integer and $\geq 2$. Then, for a small
$\epsilon>0$, let us define:
\begin{equation}
\label{phi}
\Phi^i_{\epsilon}(t,x,y)=u_{i}(t,x)-(1-\lambda)w_{i}(t,y)-\frac{1}{2\epsilon}|x-y|^{2\gamma}
-\theta(|x-\overline{x}|^{2\gamma + 2}+|y-\overline{x}|^{2\gamma +
2})-\beta (t-\overline{t})^2 - \frac{\lambda}{t}.
\end{equation}
By the growth assumption on $u_i$ and $w_i$, there exists a
$(t_{0},x_{0},y_{0},i_0)\in (0,T]\times B_R \times B_R \times {\cal
I}$, such that:
$$\Phi^{i_0}_{\epsilon}(t_{0},x_{0},y_{0})=\max\limits_{(t,x,y,i)}\Phi^i_{\epsilon}(t,x,y).$$
On the other hand, from
$2\Phi^{i_0}_{\epsilon}(t_{0},x_{0},y_{0})\geq
\Phi^{i_0}_{\epsilon}(t_0,x_0,x_0)+\Phi^{i_0}_{\epsilon}(t_0,y_0,y_0)$,
we have
\begin{equation}
\frac{1}{2\epsilon}|x_0 -y_0|^{2\gamma} \leq
(u_{i_0}(t_0,x_0)-u_{i_0}(t_0,y_0))+(1-\lambda)(w_{i_0}(t_0,x_0)-w_{i_0}(t_0,y_0)),
\end{equation}
and consequently $\frac{1}{2\epsilon}|x_0 -y_0|^{2\gamma}$ is
bounded, and as $\epsilon\rightarrow 0$, $|x_0 -y_0|\rightarrow 0$.
Since $u_{i_0}$
and $w_{i_0}$ are uniformly continuous on $[0,T]\times \overline{B}_R$, then $\frac{1}{2\epsilon}|x_0 -y_0|^{2\gamma}\rightarrow 0$ as $\epsilon\rightarrow 0.$\\
Since
 $$(1-\lambda)u_{\overline{i}}(\overline{t},\overline{x})-w_{\overline{i}}(\overline{t},\overline{x})  - \frac{\lambda}{\overline{t}}\leq
 \Phi^{i_0}_{\epsilon}(t_{0},x_{0},y_{0})\leq (1-\lambda)u_{i_0}(t_{0},x_{0})-w_{i_0}(t_{_0},y_{0})- \frac{\lambda}{t_0},$$
it follow as $\lambda\rightarrow 0$ and the continuity of $u$ and
$w$ that, up to a subsequence,
 \begin{equation}\label{subsequence}
 (t_0,x_0,y_0,i_0)\rightarrow (\overline{t},\overline{x},\overline{x},\overline{i}).
 \end{equation}

Next let us show that $t_0 <T.$ Actually if $t_0 =T$ then,
$$
\Phi^{\overline{i}}_{\epsilon}(\overline{t},\overline{x},\overline{x})\leq
\Phi^{i_0}_{\epsilon}(T,x_{0},y_{0}),$$ and,
$$
u_{\overline{i}}(\overline{t},\overline{x})-(1-\lambda)w_{\overline{i}}(\overline{t},\overline{x})-
\frac{\lambda}{\overline{t}}\leq  -\beta (T-\overline{t})^2 -
\frac{\lambda}{T},
$$
since $u_{i_0}(T,x_0)=w_{i_0}(T,y_0)=0.$ Then thanks to
(\ref{comp-equ}) we have,
$$
\begin{array}{ll}
\eta &\leq -\lambda w_{\overline{i}}(\overline{t},\overline{x})+\beta T +\frac{\lambda}{\overline{t}}\\
\eta &< \frac{3}{4}\eta .
\end{array}
$$
which yields a contradiction and we have $t_0 \in (0,T)$. We now
claim that:
\begin{equation}
\label{visco-comp1} u_{i_0}(t_0,x_0)- \left(\max\limits_{j\in{\cal
I}^{-i_0}}\{-e^{t_0}g_{i_0 j}(t_0,x_0)+u_j(t_0,x_0)\} \vee
(-e^{t_0}F_{i_0}(t_0,x_0))\right)> 0.
\end{equation}
Indeed if
$$u_{i_0}(t_0,x_0)- \left(\max\limits_{j\in{\cal I}^{-i_0}}\{-e^{t_0}g_{i_0 j}(t_0,x_0)+u_j(t_0,x_0)\}\vee(- e^{t_0}F_{i_0}(t_0,x_0))\right) \leq
0.$$

\noindent$\bf Case 1:$ $\max\limits_{j\in{\cal
I}^{-i_0}}\{-e^{t_0}g_{i_0 j}(t_0,x_0)+u_j(t_0,x_0)\}\vee(-
e^{t_0}F_{i_0}(t_0,x_0))=- e^{t_0}F_{i_0}(t_0,x_0).$ \\ \indent Then
$$u_{i_0}(t_0,x_0)\leq - e^{t_0}F_{i_0}(t_0,x_0).$$ From the supersolution property of $w_{i_0}(t_0,y_0)$, we have
$$w_{i_0}(t_0,y_0)\geq \left(\max\limits_{j\in{\cal I}^{-i_0}}\{-e^{t_0}g_{i_0 j}(t_0,y_0)+w_j(t_0,y_0)\}\vee(- e^{t_0}F_{i_0}(t_0,y_0))\right),$$ then
$$ w_{i_0}(t_0,y_0)\geq - e^{t_0}F_{i_0}(t_0,y_0).$$
It follows that:
$$u_{i_0}(t_0,x_0)-w_{i_0}(t_0,y_0)\leq -e^{t_0}( F_{i_0}(t_0,x_0)-F_{i_0}(t_0,y_0)).$$
But this contradicts the definition of (\ref{comp-equ}), since $F$,
$u$, $w$ is uniformly continuous on $[0,T]\times \overline{B}_R$ and
sending $\lambda\rightarrow0$ and the claim (\ref{visco-comp1})
holds.\\

\noindent$\bf Case 2:$ $\max\limits_{j\in{\cal
I}^{-i_0}}\{-e^{t_0}g_{i_0 j}(t_0,x_0)+u_j(t_0,x_0)\}\vee(-
e^{t_0}F_{i_0}(t_0,x_0))=\max\limits_{j\in{\cal
I}^{-i_0}}\{-e^{t_0}g_{i_0 j}(t_0,x_0)+u_j(t_0,x_0)\}.$
\\ \indent Then there exists $k \in {\cal I}^{-i_0}$ such that:
$$u_{i_0}(t_0,x_0) \leq -e^{t_0}g_{i_0
k}(t_0,x_0)+u_k(t_0,x_0).$$ From the supersolution property of
$w_{i_0}(t_0,y_0)$, we have $$w_{i_0}(t_0,y_0)\geq
\left(\max\limits_{j\in{\cal I}^{-i_0}}\{-e^{t_0}g_{i_0
j}(t_0,y_0)+w_j(t_0,y_0)\}\vee(- e^{t_0}F_{i_0}(t_0,y_0))\right),$$
then
$$ w_{i_0}(t_0,y_0)\geq -e^{t_0}g_{i_0 k}(t_0,y_0)+w_k(t_0,y_0).$$
It follows that:
$$u_{i_0}(t_0,x_0)- (1-\lambda)w_{i_0}(t_0,y_0) -(u_{k}(t_0,x_0)-(1-\lambda)w_{k}(t_0,y_0))\leq (1-\lambda)e^{t_0}g_{i_0 k}(t_0,y_0)-e^{t_0}g_{i_0 k}(t_0,x_0).$$
Now since $g_{ij}\geq \alpha >0$, for every $i\neq j$, and taking
into account of (\ref{phi}) to obtain:
$$\Phi^{i_0}_{\epsilon}(t_{0},x_{0},y_{0})-\Phi^{k}_{\epsilon}(t_{0},x_{0},y_{0})< -\alpha \lambda e^{t_0}
+e^{t_0}g_{i_0 k}(t_0,y_0)-e^{t_0}g_{i_0 k}(t_0,x_0).$$ But this
contradicts the definition of $i_0$, since $g_{i_0 k}$ is uniformly
continuous on $[0,T]\times \overline{B}_R$ and the claim
(\ref{visco-comp1}) holds.

Next let us denote
\begin{equation}
\varphi_{\epsilon}(t,x,y)=\frac{1}{2\epsilon}|x-y|^{2\gamma}+\theta(|x-\overline{x}|^{2\gamma
+ 2}+|y-\overline{x}|^{2\gamma + 2}) +\beta (t-\overline{t})^2 +
\frac{\lambda}{t}.
\end{equation}
Then we have: \be \left\{
\begin{array}{lllll}\label{derive}
D_{t}\varphi_{\epsilon}(t,x,y)=2\beta(t-\overline{t})- \frac{\lambda}{t^2},\\
D_{x}\varphi_{\epsilon}(t,x,y)=
\frac{\gamma}{\epsilon}(x-y)|x-y|^{2\gamma-2} +\theta(2\gamma + 2)
(x-\overline{x})|x-\overline{x}|^{2\gamma}, \\
D_{y}\varphi_{\epsilon}(t,x,y)=
-\frac{\gamma}{\epsilon}(x-y)|x-y|^{2\gamma-2} +
\theta(2\gamma + 2)(y-\overline{x})|y-\overline{x}|^{2\gamma},\\\\
B(t,x,y)=D_{x,y}^{2}\varphi_{\epsilon}(t,x,y)=\frac{1}{\epsilon}
\begin{pmatrix}
a_1(x,y)&-a_1(x,y) \\
-a_1(x,y)&a_1(x,y)
\end{pmatrix}+ \begin{pmatrix}
a_2(x)&0 \\
0&a_2(y)
\end{pmatrix} \\\\
\mbox{ with } a_1(x,y)=\gamma|x-y|^{2\gamma-2}I+\gamma(2\gamma -2)(x-y)(x-y)^* |x-y|^{2\gamma-4} \mbox{ and }\\
a_2(x)=\theta(2\gamma + 2)|x-\overline{x}|^{2\gamma}I+2\theta
\gamma(2\gamma + 2)(x-\overline{x})(x-\overline{x})^*
|x-\overline{x}|^{2\gamma-2 }.
\end{array}
\right. \ee Taking into account (\ref{visco-comp1}) then applying
the result by Crandall et al. (Theorem 8.3, {\cite{[CIL]}) to the
function $$
u_{i_0}(t,x)-(1-\lambda)w_{i_0}(t,y)-\varphi_{\epsilon}(t,x,y) $$ at
the point $(t_0,x_0,y_0)$, for any $\epsilon_1 >0$, we can find
$c,d\in \R$ and $X,Y \in S_k$, such that:

\be \label{lemmeishii} \left\{
\begin{array}{lllll}
(c,\frac{\gamma}{\epsilon}(x_0-y_0)|x_0-y_0|^{2\gamma-2}
+\theta(2\gamma +
2)(x_0-\overline{x})|x_0-\overline{x}|^{2\gamma},X)
\in J^{2,+}(u_{i_0}(t_0,x_0)),\\
(-d,\frac{\gamma}{\epsilon}(x_0-y_0)|x_0-y_0|^{2\gamma-2}
-\theta(2\gamma + 2)(y_0-\overline{x})|y_0-\overline{x}|^{2\gamma
},Y)\in J^{2,-}
((1-\lambda)w_{i_0}(t_0,y_0)),\\
c+d=D_{t}\varphi_{\epsilon}(t_0,x_0,y_0)=2\beta(t_0-\overline{t})- \frac{\lambda}{t_0^2} \mbox{ and finally }\\
-(\frac{1}{\epsilon_1}+||B(t_0,x_0,y_0)||)I\leq
\begin{pmatrix}
X&0 \\
0&-Y
\end{pmatrix}\leq B(t_0,x_0,y_0)+\epsilon_1 B(t_0,x_0,y_0)^2.
\end{array}
\right. \ee Taking now into account (\ref{visco-comp1}), and the
definition of viscosity solution, we get:
$$\begin{array}{l}-c+u_{i_0}(t_0,x_0)-\frac{1}{2}Tr[\sigma^{*}(t_0,x_0)X\sigma(t_0,x_0)]-\langle\frac{\gamma}{\epsilon}(x_0-y_0)|x_0-y_0|^{2\gamma-2}
+\\\qquad\qquad\qquad\qquad\qquad\theta(2\gamma +
2)(x_0-\overline{x})|x_0-\overline{x}|^{2\gamma
},b(t_0,x_0)\rangle-e^{t_0}\psi_{i_0}(t_0,x_0)\leq 0 \mbox{ and
}\\d+(1-\lambda)w_{i_0}(t_0,y_0)-\frac{1}{2}Tr[\sigma^{*}(t_0,y_0)Y\sigma(t_0,y_0)]-\langle\frac{\gamma}{\epsilon}(x_0-y_0)|x_0-y_0|^{2\gamma-2}
-\\\qquad\qquad\qquad\qquad\qquad\theta(2\gamma +
2)(y_0-\overline{x})|y_0-\overline{x}|^{2\gamma},b(t_0,y_0)\rangle-(1-\lambda)e^{t_0}\psi_{i_0}(t_0,y_0)\geq
0\end{array}$$ which implies that:
\begin{equation}
\begin{array}{llll}
\label{viscder}
-c-d +u_{i_0}(t_0,x_0)-(1-\lambda)w_{i_0}(t_0,y_0)&\leq \frac{1}{2}Tr[\sigma^{*}(t_0,x_0)X\sigma(t_0,x_0)-\sigma^{*}(t_0,y_0)Y\sigma(t_0,y_0)]\\
&\qquad +
\langle\frac{\gamma}{\epsilon}(x_0-y_0)|x_0-y_0|^{2\gamma-2},b(t_0,x_0)-b(t_0,y_0)\rangle\\&\qquad+\langle
\theta(2\gamma + 2)(x_0-\overline{x})|x_0-\overline{x}|^{2\gamma
},b(t_0,x_0)\rangle\\&\qquad +\langle \theta(2\gamma +
2)(y_0-\overline{x})|y_0-\overline{x}|^{2\gamma },b(t_0,y_0)\rangle
\\&\qquad+e^{t_0}\psi_{i_0}(t_0,x_0)-(1-\lambda)e^{t_0}\psi_{i_0}(t_0,y_0).
\end{array}
\end{equation}
But from (\ref{derive}) there exist two constants $C$ and $C_1$ such
that:
$$||a_1(x_0,y_0)||\leq C|x_0 - y_0|^{2\gamma -2} \mbox{ and }(||a_2(x_0)||\vee ||a_2(y_0)||)\leq C_1 \theta.$$
As
$$B= B(t_0,x_0,y_0)= \frac{1}{\epsilon}
\begin{pmatrix}
a_1(x_{0},y_{0})&-a_1(x_{0},y_{0}) \\
-a_1(x_{0},y_{0})&a_1(x_{0},y_{0})
\end{pmatrix}+ \begin{pmatrix}
a_2(x_0)&0 \\
0&a_2(y_0)
\end{pmatrix}$$
then
$$B\leq \frac{C}{\epsilon}|x_0 - y_0|^{2\gamma -2}
\begin{pmatrix}
I&-I \\
-I&I
\end{pmatrix}+ C_1 \theta I.$$
It follows that:
\begin{equation}
B+\epsilon_1 B^2 \leq C(\frac{1}{\epsilon}|x_0 - y_0|^{2\gamma -2}+
\frac{\epsilon_1}{\epsilon^2}|x_0 - y_0|^{4\gamma
-4})\begin{pmatrix}
I&-I \\
-I&I
\end{pmatrix}+ C_1\theta I
\end{equation}
where $C$ and $C_1$ which hereafter may change from line to line.
Choosing now $\epsilon_1=\epsilon$, yields the relation
\begin{equation}
\label{ineg_matreciel} B+\epsilon_1 B^2 \leq \frac{C}{\epsilon}(|x_0
- y_0|^{2\gamma -2}+|x_0 - y_0|^{4\gamma -4})\begin{pmatrix}
I&-I \\
-I&I
\end{pmatrix}+ C_1\theta I.
\end{equation}
Now, from (\ref{regbs1}), (\ref{lemmeishii}) and
(\ref{ineg_matreciel}) we get:
$$\frac{1}{2}Tr[\sigma^{*}(t_0,x_0)X\sigma(t_0,x_0)-\sigma^{*}(t_0,y_0)
Y\sigma(t_0,y_0)]\leq \frac{C}{\epsilon}(|x_0 - y_0|^{2\gamma}+|x_0
- y_0|^{4\gamma -2}) +C_1 \theta(1+|x_0|^2+|y_0|^2).$$ Next $$
\langle\frac{\gamma}{\epsilon}(x_0-y_0)|x_0-y_0|^{2\gamma-2},b(t_0,x_0)-b(t_0,y_0)\rangle
\leq \frac{C^2}{\epsilon}|x_0 - y_0|^{2\gamma}$$ and finally,
$$ \begin{array}{l}\langle \theta(2\gamma + 2)(x_0-\overline{x})|x_0-\overline{x}|^{2\gamma},b(t_0,x_0)\rangle +
\langle \theta(2\gamma +
2)(y_0-\overline{x})|y_0-\overline{x}|^{2\gamma },b(t_0,y_0)\rangle
\\\qquad \qquad \qquad \qquad \qquad \qquad \qquad \qquad
\leq \theta C(1+|x_0||x_0-\overline{x}|^{2\gamma + 1}+|y_0||y_0-\overline{x}|^{2\gamma + 1}).\end{array}$$ So that
by plugging into (\ref{viscder}) and note that $\lambda >0$ we
obtain:
$$\begin{array}{l}-2\beta(t_0-\overline{t})+\frac{\lambda}{t_0^2}+ u_{i_0}(t_0,x_0)-(1-\lambda)w_{i_0}(t_0,y_0)\leq \frac{C}{\epsilon}(|x_0 - y_0|^{2\gamma}+|x_0 - y_0|^{4\gamma -2})
+\\\qquad\qquad C_1 \theta (1+|x_0|^2+|y_0|^2)
+\frac{C^2}{\epsilon}|x_0 - y_0|^{2\gamma}+ \theta
C(1+|x_0||x_0-\overline{x}|^{2\gamma +
1}+|y_0||y_0-\overline{x}|^{2\gamma + 1})+\\\qquad\qquad
e^{t_0}\psi_{i_0}(t_0,x_0)-(1-\lambda)e^{t_0}\psi_{i_0}(t_0,y_0).\end{array}$$
By sending $\epsilon\rightarrow0$, $\lambda\rightarrow0$, $\theta
\rightarrow0$ and taking into account of the continuity of
$\psi_{i_0}$ and $\gamma \geq 2$, we obtain $\eta\leq 0$ which is a
contradiction. The proof of Theorem \ref{uni} is now complete.
$\Box$
\medskip

As a by-product we have the following corollary: \bcor Let
$(v^1,...,v^m)$ be a viscosity solution of (\ref{sysvi1}) which
satisfies a polynomial growth condition then for $i=1,...,m$ and
$(t,x)\in [0,T]\times \R^k$, $$ v^i(t,x)= \sup_{(\delta,\xi)\in
{\cal D}^i_t}E[\integ{t}{\gamma}\psi_{u_s}(s,X^{tx}_s)ds
-\sum_{n\geq 1}
g_{u_{\tau_{n-1}}u_{\tau_n}}(\tau_{n},X^{tx}_{\tau_{n}})
\ind_{[\tau_{n}<\gamma]}-F(\gamma,X^{tx}_\gamma,u_\gamma)\ind_{[\gamma<T]}].\Box
$$\ecor
\section{Uniqueness of the viscosity solution in infinite horizon}
Let $x\in \R^k$ and let $X^{x}$ be the solution of the following
standard SDE:
\begin{equation}\label{sde1}
dX_t^x=b(X^x_t)dt+\sigma(X^x_t)dB_t, \quad X^x_0=x
\end{equation}where the functions $b$ and $\sigma$ are the ones of
$\bf H4$. These properties of $\sigma$ and $b$ imply in particular
that $X^{x}$ solution of the standard SDE (\ref{sde1}) exists and is
unique in $\R^k $. The operator $\cal A$  defined in
(\ref{generateur2}) is the infinitesimal generator associated with
$X^{x}$.

 \beth
There are deterministic functions $v^1,...,v^m$ $:\R^k\rightarrow
\R$ such that:
$$\forall x\in \R^k,
Y_0^{i,x}=v^i(x), \,\,i=1,...,m.$$ Moreover the functions
$(v^1,...,v^m):\R^k\rightarrow \R$ are continuous, solution in
viscosity sense of the system of variational inequalities with
inter-connected obstacles (\ref{sysvi2}) and of polynomial
growth.\eeth $Proof$: The continuity of functions $v^1,...,v^m$
follows from the dynamic programming principle and is
proved in \cite{[E]}.$\Box$\\
 We are going now to address the
question of uniqueness of the viscosity solution of the system
(\ref{sysvi2}). We have the following:

\beth \label{uni2}The solution in viscosity sense of the system of
variational inequalities with inter-connected obstacles
(\ref{sysvi2}) is unique in the space of continuous functions on
$R^k$ which satisfy a polynomial growth condition, i.e., in the
space
$$\begin{array}{l}{\cal C}:=\{\varphi:  \R^k\rightarrow
\R, \mbox{ continuous and for any }\\\qquad \qquad\qquad x, \,
|\varphi(x)|\leq C(1+|x|^\mu) \mbox{ for some constants } C\quad
\mbox{and}\quad \mu\}.\end{array}$$ \eeth {\it Proof} :  We will
show by contradiction that if $u_1,...,u_m$ and $w_1,...,w_m$ are a
subsolution and a supersolution respectively for (\ref{sysvi2}) then
for any $i=1,...,m$, $u_i\leq w_i$. Therefore if we have two
solutions of (\ref{sysvi2}) then they are obviously equal. Actually
for some $R>0$ suppose there exists $(x_0,i_0)\in B_R\times {\cal
I}$ $(B_R := \{x\in \R^k; |x|<R\})$ such that:
\begin{equation}
\label{comp-equ2}
\max\limits_{(x,i)}(u_i(x)-w_i(x))=u_{i_0}(x_0)-w_{i_0}(x_0)=\eta>0.
\end{equation}
 Then, for a small
$\epsilon>0$, and $\theta,\lambda \in(0,1)$ small enough, let us
define:
\begin{equation}
\label{phi2}
\Phi^i_{\epsilon}(x,y)=u_{i}(x)-(1-\lambda)w_{i}(y)-\frac{1}{2\epsilon}|x-y|^{2\gamma}
-\theta(|x-x_0|^{2\gamma +2}+|y-x_0|^{2\gamma+2}).
\end{equation}
By the polynomial growth assumption on $u_i$ and $w_i$, there exists
a $(x_{\epsilon},y_{\epsilon},i_\epsilon)\in  B_R \times B_R \times
{\cal I}$, such that:
$$\Phi^{i_\epsilon}_{\epsilon}(x_{\epsilon},y_{\epsilon})=\max\limits_{(x,y,i)}\Phi^i_{\epsilon}(x,y).$$
On the other hand, from
$2\Phi^{i_\epsilon}_{\epsilon}(x_{\epsilon},y_{\epsilon})\geq
\Phi^{i_\epsilon}_{\epsilon}(x_\epsilon,x_\epsilon)+\Phi^{i_\epsilon}_{\epsilon}(y_\epsilon,y_\epsilon)$,
we have
\begin{equation}
\begin{array}{ll}
\frac{1}{2\epsilon}|x_\epsilon -y_\epsilon|^{2\gamma} &\leq (u_{i_\epsilon}(x_\epsilon)-u_{i_\epsilon}(y_\epsilon))+(1-\lambda)(w_{i_\epsilon}(x_\epsilon)-w_{i_\epsilon}(y_\epsilon))\\
&\leq \sum \limits_{i\in{\cal
I}}|u_{i}(x_\epsilon)-u_{i}(y_\epsilon)|+(1-\lambda)\sum
\limits_{i\in{\cal I}}|w_{i}(x_\epsilon)-w_{i}(y_\epsilon)|
\end{array}
\end{equation}

and consequently $\frac{1}{2\epsilon}|x_\epsilon
-y_\epsilon|^{2\gamma}$ is bounded, and as $\epsilon\rightarrow 0$,
$|x_\epsilon -y_\epsilon|\rightarrow 0$. Since $u_{i}$
and $w_{i}$ are uniformly continuous on $ B_R$, then $\frac{1}{2\epsilon}|x_\epsilon -y_\epsilon|^{2\gamma}\rightarrow 0$ as $\epsilon\rightarrow 0.$\\
Since
 $$u_{i_0}(x_0)-(1-\lambda)w_{i_0}(x_0) \leq \Phi^{i_\epsilon}_{\epsilon}(x_{\epsilon},y_{\epsilon})\leq u_{i_\epsilon}(x_\epsilon)-(1-\lambda)w_{i_\epsilon}(y_\epsilon),$$
 it follow as $\lambda \rightarrow 0$ and the continuity of $u_i$ and $w_i$ that, up to a subsequence,
 \begin{equation}\label{subsequence2}
 (x_\epsilon,y_\epsilon,i_\epsilon)\rightarrow (x_0,x_0,i_0).
 \end{equation}

 We now claim that:
\begin{equation}
\label{visco-comp2} u_{i_\epsilon}(x_\epsilon)-
\left(\max\limits_{j\in{\cal I}^{-i_\epsilon}}\{-g_{i_{\epsilon}
j}(x_\epsilon)+u_j(x_\epsilon)\}\vee(-F_{i_\epsilon}(x_\epsilon))\right)
> 0.
\end{equation}
Indeed if
$$u_{i_\epsilon}(x_\epsilon)- \left(\max\limits_{j\in{\cal I}^{-i_\epsilon}}\{-g_{i_{\epsilon}
j}(x_\epsilon)+u_j(x_\epsilon)\}\vee(-F_{i_\epsilon}(x_\epsilon))\right)
\leq 0.$$ \noindent$\bf Case 1:$ $\left(\max\limits_{j\in{\cal
I}^{-i_\epsilon}}\{-g_{i_{\epsilon}
j}(x_\epsilon)+u_j(x_\epsilon)\}\vee(-F_{i_\epsilon}(x_\epsilon))\right)
=-F_{i_\epsilon}(x_\epsilon).$

Then $$u_{i_\epsilon}(x_\epsilon)\leq -F_{i_\epsilon}(x_\epsilon).$$
 From the supersolution property of $w_{i_\epsilon}(y_\epsilon)$, we have
$$w_{i_\epsilon}(y_\epsilon)\geq \left(\max\limits_{j\in{\cal I}^{-i_\epsilon}}\{-g_{i_\epsilon j}(y_\epsilon)+w_j(y_\epsilon)\}\vee(-F_{i_\epsilon}
(y_\epsilon))\right),$$ then
$$ w_{i_\epsilon}(y_\epsilon)\geq -F_{i_\epsilon}(y_\epsilon).$$
It follows that:
$$u_{i_\epsilon}(x_\epsilon)-w_{i_\epsilon}(y_\epsilon)\leq -( F_{i_\epsilon}(x_\epsilon)-F_{i_\epsilon}(y_\epsilon)).$$
But this contradicts the definition of (\ref{comp-equ2}), since $F$,
$u$, $w$ is uniformly continuous on $\overline{B}_R$ and sending
$\lambda\rightarrow0$ and the claim (\ref{visco-comp2})
holds.\\

\noindent$\bf Case 2:$ $\left(\max\limits_{j\in{\cal
I}^{-i_\epsilon}}\{-g_{i_{\epsilon}
j}(x_\epsilon)+u_j(x_\epsilon)\}\vee(-F_{i_\epsilon}(x_\epsilon))\right)
=\max\limits_{j\in{\cal I}^{-i_\epsilon}}\{-g_{i_{\epsilon}
j}(x_\epsilon)+u_j(x_\epsilon)\}.$

 then there exists $k \in {\cal I}^{-i_\epsilon}$ such
that:
$$u_{i_\epsilon}(x_\epsilon) \leq -g_{i_{\epsilon}
k}(x_\epsilon)+u_k(x_\epsilon).$$ From the supersolution property of
$w_{i_\epsilon}(y_\epsilon)$, we have
$$w_{i_\epsilon}(y_\epsilon)\geq \left(\max\limits_{j\in{\cal I}^{-i_\epsilon}}\{-g_{i_\epsilon j}(y_\epsilon)+w_j(y\epsilon)\}\vee(-F_{i_\epsilon}
(y_\epsilon))\right),$$ then
$$ w_{i_\epsilon}(y_\epsilon)\geq -g_{i_{\epsilon} k}(y_\epsilon)+w_k(y_\epsilon).$$
It follows that:
$$u_{i_\epsilon}(x_\epsilon)- (1-\lambda)w_{i_\epsilon}(y_\epsilon) -(u_{k}(x_\epsilon)-(1-\lambda)w_{k}(y_\epsilon))\leq (1-\lambda)g_{i_{\epsilon} k}(y_\epsilon)-g_{i_{\epsilon} k}(x_\epsilon).$$
Now since $g_{ij}\geq \alpha >0$, for every $i\neq j$, and taking
into account of (\ref{phi2}) to obtain:
$$
\begin{array}{ll}
\Phi^{i_\epsilon}_{\epsilon}(x_{\epsilon},y_{\epsilon})-\Phi^{k}_{\epsilon}(x_{\epsilon},y_{\epsilon})&<
-\alpha \lambda
+g_{i_{\epsilon} k}(y_\epsilon)-g_{i_{\epsilon} k}(x_\epsilon)\\
\end{array}
$$
But this contradicts the definition of $i_\epsilon$, since
$g_{i_{\epsilon} k}$ is uniformly continuous on $\overline{B}_R$ and
the claim (\ref{visco-comp2}) holds.

Next let us denote
\begin{equation}
\varphi_{\epsilon}(x,y)=\frac{1}{2\epsilon}|x-y|^{2\gamma}
+\theta(|x-x_0|^{2\gamma +2}+|y-x_0|^{2\gamma+2}).
\end{equation}
Then we have: \be \left\{
\begin{array}{lll}\label{derive2}
D_{x}\varphi_{\epsilon}(t,x,y)=
\frac{\gamma}{\epsilon}(x-y)|x-y|^{2\gamma-2} +\theta(2\gamma + 2)
(x-x_0)|x-x_0|^{2\gamma}, \\
D_{y}\varphi_{\epsilon}(t,x,y)=
-\frac{\gamma}{\epsilon}(x-y)|x-y|^{2\gamma-2} +
\theta(2\gamma + 2)(y-y_0)|y-y_0|^{2\gamma},\\\\
B(t,x,y)=D_{x,y}^{2}\varphi_{\epsilon}(t,x,y)=\frac{1}{\epsilon}
\begin{pmatrix}
a_1(x,y)&-a_1(x,y) \\
-a_1(x,y)&a_1(x,y)
\end{pmatrix}+ \begin{pmatrix}
a_2(x)&0 \\
0&a_2(y)
\end{pmatrix} \\\\
\mbox{ with } a_1(x,y)=\gamma|x-y|^{2\gamma-2}I+\gamma(2\gamma -2)(x-y)(x-y)^* |x-y|^{2\gamma-4} \mbox{ and }\\
a_2(x)=\theta(2\gamma + 2)|x-x_0|^{2\gamma}I+2\theta \gamma(2\gamma
+ 2)(x-x_0)(x-x_0)^* |x-x_0|^{2\gamma-2 }
\end{array}
\right. \ee where $I$ stands for the identity matrix of dimension
$k$. Taking into account (\ref{visco-comp2}) then applying the
result by Crandall et al. (Theorem 3.2, {\cite{[CIL]}) to the
function $$ u_{i}(x)-(1-\lambda)w_{i}(y)-\varphi_{\epsilon}(x,y) $$
at the point $(x_\epsilon,y_\epsilon)$, for any $\epsilon_1 >0$, we
can find
 $X,Y \in S_k$, such that:

\be \label{lemmeishii2} \left\{
\begin{array}{lllll}
(\frac{\gamma}{\epsilon}(x_\epsilon-y_\epsilon)|x_\epsilon-y_\epsilon|^{2\gamma-2}
+\theta(2\gamma + 2) (x_\epsilon-x_0)|x_\epsilon-x_0|^{2\gamma},X)
\in J^{2,+}(u_{i_\epsilon}(x_\epsilon)),\\
(\frac{\gamma}{\epsilon}(x_\epsilon-y_\epsilon)|x_\epsilon-y_\epsilon|^{2\gamma-2}
- \theta(2\gamma +
2)(y_\epsilon-y_0)|y_\epsilon-y_0|^{2\gamma},Y)\in J^{2,-}
((1-\lambda)w_{i_\epsilon}(y_\epsilon)),\\
-(\frac{1}{\epsilon_1}+||B(x_\epsilon,y_\epsilon)||)\begin{pmatrix}
I&0 \\
0&I
\end{pmatrix}\leq
\begin{pmatrix}
X&0 \\
0&-Y
\end{pmatrix}\leq B(x_\epsilon,y_\epsilon)+\epsilon_1 B(x_\epsilon,y_\epsilon)^2.
\end{array}
\right. \ee Taking now into account (\ref{visco-comp2}), and the
definition of viscosity solution, we get:
$$\begin{array}{l}ru_{i_\epsilon}(x_\epsilon)-\frac{1}{2}Tr[\sigma^{*}(x_\epsilon)X\sigma(x_\epsilon)]-\langle\frac{\gamma}{\epsilon}
(x_\epsilon-y_\epsilon)|x_\epsilon-y_\epsilon|^{2\gamma-2}
\\\qquad\qquad\qquad\qquad\qquad +\theta(2\gamma + 2)
(x_\epsilon-x_0)|x_\epsilon-x_0|^{2\gamma},b(x_\epsilon)\rangle-\psi_{i_\epsilon}(x_\epsilon)\leq
0 \mbox{ and
}\\r(1-\lambda)w_{i_\epsilon}(y_\epsilon)-\frac{1}{2}Tr[\sigma^{*}(y_\epsilon)Y\sigma(y_\epsilon)]-\langle
\frac{\gamma}{\epsilon}
(x_\epsilon-y_\epsilon)|x_\epsilon-y_\epsilon|^{2\gamma-2}
\\\qquad\qquad\qquad\qquad\qquad -\theta(2\gamma + 2)
(y_\epsilon-x_0)|y_\epsilon-x_0|^{2\gamma},b(y_\epsilon)\rangle-(1-\lambda)\psi_{i_\epsilon}(y_\epsilon)\geq
0\end{array}$$ which implies that:
\begin{equation}
\begin{array}{llll}
\label{viscder2}
&ru_{i_\epsilon}(x_\epsilon)-r(1-\lambda)w_{i_\epsilon}(y_\epsilon)\leq \frac{1}{2}Tr[\sigma^{*}(x_\epsilon)X\sigma(x_\epsilon)-\sigma^{*}(y_\epsilon)Y\sigma(y_\epsilon)]\\
&\qquad + \langle\frac{\gamma}{\epsilon}
(x_\epsilon-y_\epsilon)|x_\epsilon-y_\epsilon|^{2\gamma-2},b(x_\epsilon)-b(y_\epsilon)\rangle\\&\qquad+\langle
\theta(2\gamma + 2)
(x_\epsilon-x_0)|x_\epsilon-x_0|^{2\gamma},b(x_\epsilon)\rangle
+\langle \theta(2\gamma + 2)
(y_\epsilon-x_0)|y_\epsilon-x_0|^{2\gamma},b(y_\epsilon)\rangle
\\&\qquad+\psi_{i_\epsilon}(x_\epsilon)-(1-\lambda)\psi_{i_\epsilon}(y_\epsilon).
\end{array}
\end{equation}
But from (\ref{derive2}) there exist two constants $C$ and $C_1$
such that:
$$||a_1(x_\epsilon,y_\epsilon)||\leq C|x_\epsilon - y_\epsilon|^{2\gamma -2} \mbox{ and }
(||a_2(x_\epsilon)||\vee ||a_2(y_\epsilon)||)\leq C_1\theta .$$ As
$$B= B(x_\epsilon,y_\epsilon)= \frac{1}{\epsilon}
\begin{pmatrix}
a_1(x_\epsilon,y_\epsilon)&-a_1(x_\epsilon,y_\epsilon) \\
-a_1(x_\epsilon,y_\epsilon)&a_1(x_\epsilon,y_\epsilon)
\end{pmatrix}+ \begin{pmatrix}
a_2(x_\epsilon)&0 \\
0&a_2(y_\epsilon)
\end{pmatrix}$$
then
$$B\leq \frac{C}{\epsilon}|x_\epsilon - y_\epsilon|^{2\gamma -2}
\begin{pmatrix}
I&-I \\
-I&I
\end{pmatrix}+ C_1\theta \begin{pmatrix}
I&0 \\
0&I
\end{pmatrix}.$$
It follows that:
\begin{equation}
B+\epsilon_1 B^2 \leq C(\frac{1}{\epsilon}|x_\epsilon -
y_\epsilon|^{2\gamma -2}+ \frac{\epsilon_1}{\epsilon^2}|x_\epsilon -
y_\epsilon|^{4\gamma -4})\begin{pmatrix}
I&-I \\
-I&I
\end{pmatrix}+ C_1\theta \begin{pmatrix}
I&0 \\
0&I
\end{pmatrix}
\end{equation}
where $C$ and $C_1$ which hereafter may change from line to line.
Choosing now $\epsilon_1=\epsilon$, yields the relation
\begin{equation}
\label{ineg_matreciel2} B+\epsilon_1 B^2 \leq
\frac{C}{\epsilon}(|x_\epsilon - y_\epsilon|^{2\gamma
-2}+|x_\epsilon - y_\epsilon|^{4\gamma -4})\begin{pmatrix}
I&-I \\
-I&I
\end{pmatrix}+ C_1\theta \begin{pmatrix}
I&0 \\
0&I
\end{pmatrix}.
\end{equation}
Now, from $\bf H4$, (\ref{lemmeishii2}) and (\ref{ineg_matreciel2})
we get:
$$\frac{1}{2}Tr[\sigma^{*}(x_\epsilon)X\sigma(x_\epsilon)-\sigma^{*}(y_\epsilon)
Y\sigma(y_\epsilon)]\leq \frac{C}{\epsilon}(|x_\epsilon -
y_\epsilon|^{2\gamma}+|x_\epsilon - y_\epsilon|^{4\gamma -2}) +C_1
\theta(1+|x_\epsilon|^2+|y_\epsilon|^2).$$ Next $$
\langle\frac{\gamma}{\epsilon}
(x_\epsilon-y_\epsilon)|x_\epsilon-y_\epsilon|^{2\gamma-2},b(x_\epsilon)-b(y_\epsilon)\rangle
\leq \frac{C^2}{\epsilon}|x_\epsilon - y_\epsilon|^{2\gamma}$$ and
finally,
$$\langle
\theta(2\gamma + 2)
(x_\epsilon-x_0)|x_\epsilon-x_0|^{2\gamma},b(x_\epsilon)\rangle \leq
\theta C(1+|x_\epsilon|)|x_\epsilon -x_0|^{2\gamma+1}
$$
$$
\langle \theta(2\gamma + 2)
(y_\epsilon-x_0)|y_\epsilon-x_0|^{2\gamma},b(y_\epsilon)\rangle \leq
\theta C(1+|y_\epsilon|)|y_\epsilon -x_0|^{2\gamma+1}.$$ So that by
plugging into (\ref{viscder2}) we obtain:
$$\begin{array}{l}ru_{i_\epsilon}(x_\epsilon)-r(1-\lambda)w_{i_\epsilon}(y_\epsilon) \leq \frac{C}{\epsilon}(|x_\epsilon -
y_\epsilon|^{2\gamma}+|x_\epsilon - y_\epsilon|^{4\gamma -2}) +C_1
\theta(1+|x_\epsilon|^2+|y_\epsilon|^2)+
\\\qquad \qquad \frac{C^2}{\epsilon}|x_\epsilon - y_\epsilon|^{2\gamma}+ \theta C(1+|x_\epsilon|)|x_\epsilon -x_0|^{2\gamma+1}+\theta
C(1+|y_\epsilon|)|y_\epsilon -x_0|^{2\gamma+1}+\\\qquad \qquad
\psi_{i_\epsilon}(x_\epsilon)-(1-\lambda)\psi_{i_\epsilon}(y_\epsilon).\end{array}$$
By sending $\epsilon\rightarrow0$, $\lambda\rightarrow0$, $\theta
\rightarrow0$ and taking into account of the continuity of
$\psi_{i_\epsilon}$, we obtain $u_{i_0}(x_0)-w_{i_0}(x_0)<0$ which
is a contradiction. The proof of Theorem \ref{uni2} is now
complete.$\Box$

As a by-product we have the following result: \bcor Let
$(v^1,...,v^m)$ be a viscosity solution of (\ref{sysvi2}) which
satisfies a polynomial growth condition. Then for $i=1,...,m$ and
$(t,x)\in \R^k$, $$ v^i(x)= \sup_{(\delta,\xi)\in {\cal
D}^i_0}E[\integ{0}{\gamma}e^{-rs}\psi_{u_s}(X^{x}_s)ds -\sum_{n\geq
1}
e^{-r\tau_n}g_{u_{\tau_{n-1}}u_{\tau_n}}(X^{x}_{\tau_{n}})\ind_{[\tau_n<\gamma]}
-e^{-r\gamma}F(X^x_\gamma,u_\gamma)].\Box
$$ \ecor


\end{document}